\documentclass[12pt]{amsart}

\usepackage{times,amsfonts,amsmath,amstext,amsbsy,amssymb,
  amsopn,amsthm,upref,eucal, amscd, graphicx}
\usepackage[pdftex]{color}
\usepackage[T1]{fontenc}

\newtheorem{theorem}{Theorem}[section]

\numberwithin{equation}{section}

\theoremstyle{definition}

\newtheorem{remark}[theorem]{Remark}

\begin{document}
\title[Non-uniformly hyperbolic horseshoes in the standard family]{Non-uniformly hyperbolic horseshoes in the standard family}
\author{Carlos Matheus}
\address{Carlos Matheus: Universit\'e Paris 13, Sorbonne Paris Cit\'e, LAGA, CNRS (UMR 7539), F-93439, Villetaneuse, France}
\email{matheus@math.univ-paris13.fr.}
\author{Carlos Gustavo Moreira}
\address{Carlos Gustavo Moreira: IMPA, Estrada D. Castorina, 110, CEP 22460-320, Rio de Janeiro, RJ, Brazil}
\email{gugu@impa.br.}
\author{Jacob Palis}
\address{Jacob Palis: IMPA, Estrada D. Castorina, 110, CEP 22460-320, Rio de Janeiro, RJ, Brazil}
\email{jpalis@impa.br.}
\date{\today}
\begin{abstract}
We show that the non-uniformly hyperbolic horseshoes of Palis and Yoccoz occur in the standard family of area-preserving diffeomorphisms of the two-torus. 
\end{abstract}
\maketitle


\section{Introduction}\label{intro}

In their \emph{tour-de-force} work about the dynamics of surface diffeomorphisms, Palis and Yoccoz \cite{PY} proved that the so-called \emph{non-uniformly hyperbolic horseshoes} are very frequent in the generic unfolding of a first heteroclinic tangency associated to periodic orbits in a horseshoe with Hausdorff dimension slightly bigger than one. 

In the same article, Palis and Yoccoz gave an \emph{ad-hoc} example of $1$-parameter family of diffeomorphisms of the two-sphere fitting the setting of their main results, and, thus, exhibiting non-uniformly hyperbolic horseshoes: see page 3 (and, in particular, Figure 1) of \cite{PY}. 

In this note, we show that the \emph{standard family} $f_k:\mathbb{T}^2\to \mathbb{T}^2$, $k\in\mathbb{R}$, 
$$f_k(x,y) := (-y+2x+k\sin(2\pi x), x)$$ 
of area-preserving diffeomorphisms of the two-torus $\mathbb{T}^2=\mathbb{R}^2/\mathbb{Z}^2$ displays non-uniformly hyperbolic horseshoes. 

More precisely, our main theorem is: 

\begin{theorem}\label{t.MMP} There exists $k_0 > 0$ such that, for all $|k|>k_0$, the subset of parameters $r\in\mathbb{R}$ such that $|r-k|<4/k^{1/3}$ and $f_r$ exhibits a non-uniformly hyperbolic horseshoe (in the sense of Palis--Yoccoz \cite{PY}) has positive Lebesgue measure. 
\end{theorem}

The remainder of this text is divided into three sections: in Section \ref{s.PYsetting}, we briefly recall the context of Palis--Yoccoz work \cite{PY}; in Section \ref{s.Duarte}, we revisit some elements of Duarte's construction \cite{Duarte} of tangencies associated to certain (uniformly hyperbolic) horseshoes of $f_k$; finally, we establish Theorem \ref{t.MMP} in Section \ref{s.MMP} by modifying Duarte's constructions (from Section \ref{s.Duarte}) in order to apply Palis-Yoccoz results (from Section \ref{s.PYsetting}). 

\section{Non-uniformly hyperbolic horseshoes}\label{s.PYsetting} 

Suppose that $F$ is a smooth diffeomorphism of a compact surface $M$ displaying a first heteroclinic tangency associated to periodic points of a horseshoe $K$, that is:
\begin{itemize}
\item $p_s, p_u\in K$ belong to distinct periodic orbits of $F$; 
\item $W^s(p_s)$ and $W^u(p_u)$ have a quadratic tangency at a point $q\in M\setminus K$; 
\item for some neighborhoods $U$ of $K$ and $V$ of the orbit $\mathcal{O}(q)$ such that the maximal invariant set of $U\cup V$ is precisely $K\cup\mathcal{O}(q)$. 
\end{itemize}

Assume that $K$ is \emph{slightly thick} in the sense that its stable and unstable dimensions $d^s$ and $d^u$ satisfy $d_s+d_u>1$ and 
$$(d_s+d_u)^2 + \max(d_s+d_u)^2 < d_s+d_u+\max(d_s, d_u)$$

\begin{remark}\label{r.conservative-horseshoe} Since the stable and unstable dimensions of a horseshoe of an \emph{area-preserving} diffeomorphism $F$ always coincide, a slightly thick horseshoe $K$ of an area-preserving diffeormorphism $F$ has stable and unstable dimensions 
$$0.5< d_s = d_u < 0.6$$
\end{remark}

In this setting, the results proved by Palis and Yoccoz \cite{PY} imply the following statement: 

\begin{theorem}[Palis--Yoccoz]\label{t.PY} Given a $1$-parameter family $(F_t)_{|t|<t_0}$ with $F_0 = F$ and generically unfolding the heteroclinic tangency at $q$, the subset of parameters $t\in (-t_0, t_0)$ such that $F_t$ has a non-uniformly hyperbolic horseshoe\footnote{We are not going to recall the definition of non-uniformly hyperbolic horseshoes here: instead, we refer to the original article \cite{PY} for the details.} has positive Lebesgue measure.  
\end{theorem} 

\section{Horseshoes and tangencies in the standard family}\label{s.Duarte}

The standard family $f_k$ generically unfolds tangencies associated to very thick horseshoes $\Lambda_k$: this phenomenon was studied in details by Duarte \cite{Duarte} during his proof of the almost denseness of elliptic islands of $f_k$ for large generic parameters $k$.

In the sequel, we review some facts from Duarte's article about $\Lambda_k$ and its tangencies (for later use in the proof of our Theorem \ref{t.MMP}). 

For technical reasons, it is convenient to work with the standard family $f_k$ and their \emph{singular} perturbations 
$$g_k(x,y)=(-y+2x+k\sin(2\pi x)+\rho_k(x), x),$$ 
where $\rho_k$ is defined in Section 4 of \cite{Duarte}. Here, it is worth to recall that the key features of $\rho_k$ are: 
\begin{itemize}
\item $\rho_k$ has \emph{poles} at the critical points $\nu_{\pm}=\pm 1/4+O(1/k)$ of the function $2x+k\sin(2\pi x)$; 
\item $\rho_k$ vanishes outside $|x\pm\frac{1}{4}|\leq \frac{2}{k^{1/3}}$. 
\end{itemize} 

In Section 2 of \cite{Duarte}, Duarte constructs the stable and unstable foliations $\mathcal{F}^s$ and $\mathcal{F}^u$ for $g_k$. As it turns out, $\mathcal{F}^s$, resp. $\mathcal{F}^u$, is an almost vertical, resp. horizontal, foliation in the sense that it is generated by a vector field $(\alpha^s(x,y),1)$, resp. $(1,\alpha^u(x,y))$, satisfying all properties described in Section 2 of Duarte's paper \cite{Duarte}. In particular, $\mathcal{F}^s$, resp. $\mathcal{F}^u$, describe the local stable, resp. unstable, manifolds for the standard map $f_k$ at points whose future, resp. past, orbits stay in the region $\{f_k=g_k\}$, resp. $\{f_k^{-1}=g_k^{-1}\}$. 

In Section 3 of \cite{Duarte}, Duarte analyses the projections $\pi^s$ and $\pi^u$ obtained by thinking the foliations $\mathcal{F}^s$ and $\mathcal{F}^u$ as fibrations over the singular circles $C_s=\{(x,\nu_+)\in\mathbb{T}^2\}$ and $C_u=\{(\nu_+,y)\in\mathbb{T}^2\}$. Among many things, Duarte shows that the circle map $\Psi_k:\mathbb{S}^1\to\mathbb{S}^1$ defined by 
$$(\Psi_k(x), \nu_+):=\pi^s(g_k(x,\nu_+)) \textrm{ or, equivalently, } (\nu_+,\Psi_k(y)) = \pi^u(g_k^{-1}(\nu_+,y))$$ 
is \emph{singular} expansive with small distortion. 

In Section 4 of \cite{Duarte}, Duarte considers a Cantor set 
$$K_k=\bigcap_{n\in\mathbb{N}}\Psi_k^{-1}(J_0\cup J_1)$$ 
of the circle map $\Psi_k$ associated to a Markov partition $J_0\cup J_1 \subset [-1/4, 3/4]$ with the following properties:
\begin{itemize}
\item the extremities of the intervals $J_0=[a,b]$ and  $J_1=[b', a'+1]$ satisfy $a+\frac{1}{4}, \frac{1}{4}-b, -\frac{1}{4}-a', b'-\frac{1}{4}, \in (\frac{3}{k^{1/3}}, \frac{4}{k^{1/3}})$, so that $J_0$ and $J_1$ are contained in the region $\{\rho_k=0\}$;  
\item $\Psi_k(a) = a = \Psi_k(a')$, $\Psi_k(b) = a' = \Psi_k(b')$.
\end{itemize} 
In particular, Duarte uses these features of $K_k$ to prove that 
$$\Lambda_k = (\pi^s)^{-1}(K_k)\cap (\pi^u)^{-1}(K_k)$$ 
is a horseshoe of \emph{both} $g_k$ and $f_k$. 

In Section 5 of \cite{Duarte}, Duarte studies the tangencies associated to the invariant foliations of $\Lambda_k$. More concretely, denote by $\mathcal{G}^u = (f_k)_*(\mathcal{F}^u)$ the foliation obtained by pushing the almost horizontal foliation $\mathcal{F}^u$ by the standard map $f_k$. The vector fields $(\beta^u(x,y),1)$ defining $\mathcal{G}^u$ and $(\alpha^s(x,y), 1)$ defining $\mathcal{F}^s$ coincide along two (almost horizontal) circles of tangencies $\{(x,\sigma_+(x):x\in\mathbb{S}^1\}\cup\{(x,\sigma_-(x):x\in\mathbb{S}^1\}$ (with $|\sigma_{\pm}(x)-\nu_{\pm}|\leq \frac{1}{270 k^{5/3}}$ and $|\sigma_{\pm}'(x)|\leq \frac{1}{12 k^{4/3}}$ for all $x\in\mathbb{S}^1$). The projections of 
$\Lambda_k$ along $\mathcal{F}^s$ and $\mathcal{G}^u$ on the circle of tangencies $\{(x,\sigma_+(x)):x\in\mathbb{S}^1\}$ define two Cantor sets 
$$K^s_h = \{(x,\sigma_+(x)):x\in\mathbb{S}^1\}\cap (\pi^s)^{-1}(K_k)$$
and 
$$K^u_h = \{(x,\sigma_+(x)):x\in\mathbb{S}^1\}\cap f_k((\pi^u)^{-1}(K_k))$$ 
whose intersection points $x\in K^s_h\cap K^u_h$ are points of tangencies between the invariant manifolds of $\Lambda_k$. Furthermore, it is shown in Propositions 18 and 20 of \cite{Duarte} that these tangencies are quadratic\footnote{The difference in curvatures at tangency points is $\geq 4\pi^2k - \frac{3}{k^{1/3}}$} and unfold generically\footnote{The leaves of $\mathcal{F}^s$ move with speed $\leq \frac{3}{k^{2/3}}$ and the leaves of $\mathcal{G}^u$ move with speed $\geq 1-\frac{3}{k^{2/3}}$.}\label{f.speed}. 

\section{Proof of Theorem \ref{t.MMP}}\label{s.MMP}

After these preliminaries on the works of Palis--Yoccoz and Duarte, we are ready to prove the main result of this note. 

The standard map $f_k$ has fixed points at $p_s = (0,0)\in\Lambda_k$ and $p_u = (-\frac{1}{12}+O(\frac{1}{k}), -\frac{1}{12}+O(\frac{1}{k}))\in\Lambda_k$.

The local stable leaf $\mathcal{F}^s(p_s)$ is tangent to some leaf of $\mathcal{G}^u$ at a point $q$. Since $K_k$ is $\frac{2}{k^{1/3}}$-dense in $\mathbb{S}^1$ (cf. page 394 of \cite{Duarte}), and $f_k$ sends the vertical circle $f_k^{-1}(\{(x,\sigma_+(x)):x\in\mathbb{S}^1\}) := \{(\rho_+(x), x):x\in\mathbb{S}^1\}$ into the horizontal circle $\{(x,\sigma_+(x)):x\in\mathbb{S}^1\}$ as a $C^1$-perturbation of size $\frac{1}{81 k^2}$ of a rigid rotation (cf. page 397 of \cite{Duarte}), we can find a point of $K^u_h$ in the $\frac{7}{2k^{1/3}}$-neighborhood of the tangency point $q\in \{(x,\sigma_+(x)):x\in\mathbb{S}^1\}$. 

Therefore, the fact that the tangency at $q$ unfolds generically (cf. footnote \ref{f.speed}) permits to take a parameter $|\overline{k}-k|<\frac{4}{k^{1/3}}$ such that the local stable leaf $\mathcal{F}^s(p_s)$ is tangent to the unstable manifold of some point of $\Lambda_{\overline{k}}$. 

Because the unstable manifold of the fixed point $p_u$ is dense in $\Lambda_{\overline{k}}$ (and the tangencies unfold generically), we can replace $\overline{k}$ by a parameter $|r-k|<\frac{4}{k^{1/3}}$ such that the local stable manifold $\mathcal{F}^s(p_s)$ has a quadratic tangency with the unstable manifold of $p_u$ at $q$ which is unfolded generically. 

Next, we observe that the right part of a small neighborhood of $q$ in the circle of tangencies is transversal to leaves of $\mathcal{F}^s$ to the right of $p_s$, and the left part of a small neighborhood of $q$ in the circle of tangencies is transversal to a certain (fixed) iterate of the leaves of $\mathcal{F}^u$ which are either all above or all below $p_u$. In the former, resp. latter, case, we consider a Markov partition $I_{-}\cup I_0\cup I_1$ for the singular expansive map $\Psi_r:\mathbb{S}^1\to \mathbb{S}^1$ where: 
\begin{itemize} 
\item $I_0$ has extremities $\pi^s(p_s)$ and $a\in[\frac{1}{8}, \frac{1}{8}+\frac{1}{k^{1/3}}]$; 
\item $I_1$ has extremities $b\in [\frac{15}{32}-\frac{1}{k^{1/3}}, \frac{15}{32}]$ and $c\in[\frac{19}{32}, \frac{19}{32}+\frac{1}{k^{1/3}}]$;
\item $I_-$ has extremities $\pi^u(p_u)$ and $d\in[-\frac{1}{48}, -\frac{1}{48}+\frac{1}{k^{1/3}}]$, resp. $d\in [-\frac{7}{48}-\frac{1}{k^{1/3}}, -\frac{7}{48}]$;
\item $\Psi_r(c)=\pi^u(p_u)$, $\Psi_r(b)=c=\Psi_r(d)$ and $\Psi_r(a)=d$, resp. $\Psi_r(a)=\pi^u(p_u)$, $\Psi_r(d) = \pi^s(p_s)$, $\Psi_r(c)=d$ and $\Psi_r(b) = c$. 
\end{itemize} 
This defines a Cantor set 
$$L_r:=\bigcap\limits_{n\in\mathbb{N}} \Psi_r^{-n}(I_-\cup I_0\cup I_1)$$ 
and a horseshoe 
$$\Theta_r:=(\pi^s)^{-1}(L_r)\cap(\pi^u)^{-1}(L_r)$$ 
containing $p_s$ and $p_u$. 

By definition, we can select neighborhoods $U$ of $\Theta_r$ and $V$ of the orbit $\mathcal{O}(q)$ of $q$ such that the $f_r$-maximal invariant set of $U\cup V$ is exactly $\Theta_r\cup\mathcal{O}(q)$: this happens because our choices were made so that the local stable leafs of $\Theta_r$ approach $q$ only from the right, while certain (fixed) iterates of the local unstable manifolds of $\Theta_r$ approach $q$ only from the left. 

Therefore, we can conclude Theorem \ref{t.MMP} from Palis--Yoccoz work (cf. Theorem \ref{t.PY}) once we verify that $\Theta_r$ is slightly thick. 

In view of Remark \ref{r.conservative-horseshoe}, our task is reduced to check that the stable and unstable Hausdorff dimensions of $\Theta_r$ are comprised between $0.5$ and $0.6$. In this direction, note that these Hausdorff dimensions coincide with the Hausdorff dimension $d(r)$ of $L_r$. Moreover, the distortion constant $C_1(r)$ of $\Psi_r$ is small (namely, $0\leq C_1(k)\leq \frac{9}{k^{1/3}}$, cf. page 388 of \cite{Duarte}). Hence, $d(r)$ is close to the solution $\kappa(r)$ of ``Bowen's equation'' 
$$(\textrm{length } I_-)^{\kappa(r)} + (\textrm{length } I_0)^{\kappa(r)} + (\textrm{length } I_1)^{\kappa(r)} = (\textrm{length } I)^{\kappa(r)}$$
where $I$ is the convex hull of $I_-\cup I_0\cup I_1$. Since $\textrm{length } I_-=\frac{1}{16}+O(\frac{1}{k^{1/3}})$, $\textrm{length } I_0 =\textrm{length } I_1 = \frac{1}{8}+O(\frac{1}{k^{1/3}})$, 
$$\textrm{length } I = \frac{19}{32}+\frac{1}{12}+O(\frac{1}{k^{1/3}}), \quad \textrm{resp. } \,\, \frac{19}{32}+\frac{7}{48}+O(\frac{1}{k^{1/3}})$$
and 
$$\left(1/16\right)^{0.5809\dots} + \left(1/8\right)^{0.5809\dots} + \left(1/8\right)^{0.5809\dots} = \left(65/96\right)^{0.5809\dots}, \textrm{ resp.}$$
$$\left(1/16\right)^{0.5546\dots} + \left(1/8\right)^{0.5546\dots} + \left(1/8\right)^{0.5546\dots} = \left(71/96\right)^{0.5546\dots},$$
we derive that $0.554<d(r)<0.581$. This completes the argument.

\end{document}